\documentclass{elsarticle}

\usepackage{amssymb}
\usepackage{amsmath}
\usepackage{enumerate}

\usepackage[T1]{fontenc}
\usepackage[utf8]{inputenc}

\usepackage{lmodern}
\usepackage[english,francais]{babel}

\def\sqw{\hbox{\rlap{\leavevmode\raise.3ex\hbox{$\sqcap$}}$%
\sqcup$}}

\def\cqfd{\ifmmode\sqw\else{\ifhmode\unskip\fi\nobreak\hfil
\penalty50\hskip1em\null\nobreak\hfil\sqw
\parfillskip=0pt\finalhyphendemerits=0\endgraf}\fi}

\def \dem{\noindent{\sc\small\it{ Proof.} }}

\newtheorem{theorem}{Theorem}[section]
\newtheorem{lemma}[theorem]{Lemma}
\newtheorem{proposition}[theorem]{Proposition}
\newtheorem{corollary}[theorem]{Corollary}

\newtheorem{remark}[theorem]{Remark}

\def\og{\leavevmode\raise.3ex\hbox{$\scriptscriptstyle\langle\!\langle$~}}
\def\fg{\leavevmode\raise.3ex\hbox{~$\!\scriptscriptstyle\,\rangle\!\rangle$}}

\newcommand{\N}{\mathbb{N}}
\newcommand{\Z}{\mathbb{Z}}
\newcommand{\C}{\mathbb{C}}
\def\germ #1 {\mathfrak{#1}}
\def\cal #1 {\mathcal{#1}}

\journal{Journal of Algebra}

\begin{document}

\begin{frontmatter}





 \title{A Howe-type correspondence for the dual pair $(\germ sl _2, \germ sl _n)$ in $\germ sl _{2n}$}
 \author{Guillaume Tomasini}
 \ead{tomasini@math.u-strasbg.fr}
 \address{Institut de Recherche Math\'ematique Avanc\'ee, UMR 7501\\
Universit\'e de Strasbourg et CNRS\\
7 rue Ren\'e Descartes, 67084 Strasbourg Cedex, France}


 \medskip


 \begin{abstract} 



  
In this article, we study the decomposition of weight--$\germ sl _{2n}$--modules of degree $1$ to a dual pair $(\germ sl _2,\germ sl _n)$. We show that in some generic cases we have an explicit branching rule leading to a Howe--type correspondence between simple highest weight modules. We also give a Howe--type correspondence in the non--generic case. This latter involves some (non simple) Verma modules.

 \end{abstract}
\end{frontmatter}

\selectlanguage{english}


Let $\germ g $ denote a reductive Lie algebra over $\C$. A dual pair in $\germ g $ is a pair $(\germ a , \germ b )$ of reductive subalgebras of $\germ g $ which are the commutant of each other. Given a simple $\germ g $--module $M$, one can try to solve the following branching problem: describe the restriction of $M$ to the subalgebra $\germ a + \germ b $. This problem and his group analogue have received particular attention since the late $80s$. The first result concerning such a restriction was obtained by R. Howe in \cite{Ho89a} and \cite{Ho89b}. These articles were concerned with the Lie algebra $\germ sp _{2n}$ (in fact the metaplectic group whose Lie algebra is the symplectic Lie algebra) and the so--called minimal (or Weil, or Shale--Segal--Weil, or oscillator) representation. From the infinitesimal point of view, the vector space of the representation is a polynomial algebra and the action is via differential operators. The restriction of this representation to the dual pair gives rise to a one--one correspondence between some simple representations of $\germ a $ and some simple representations of $\germ b $. The correspondence from the point of view of Lie group, and for the Weil representation, is usually called $\theta$--correspondence. In the case of Lie algebras we call such a correspondence a Howe--type correspondence or a dual pair correspondence. Other occurences of such a correspondence can be found in \cite{RS95}, \cite{HPS96}, \cite{Li99b}. All these articles deal with the minimal representation of some Lie algebra or Lie group.

The aim of this article is to prove a Howe--type correspondence for a new familly of representations of the Lie algebra $\germ sl _{2n}$, which was introduced by Benkart, Britten, and Lemire in \cite{BBL97}. The vector space of the representation is some kind of polynomial algebra and the action is via differential operators. The correspondence is completely explicit (see theorems \ref{thmAB1}, \ref{thmAB2} and \ref{thmAB3}).

In the first part of this article we give the construction of the representation and some of its properties. The second part is devoted to the description of the dual pair $(\germ sl _2, \germ sl _{n})$ of $\germ sl _{2n}$ and its action on the representation of the first part. In the last part we prove the Howe--type correspondence for this module with respect to our dual pair.

{\it Acknowledgements.--} I thank gratefully Professor H. Rubenthaler for many helpful conversations and valuable comments concerning the writing of this article.



\section{Simple weight--modules of degree $1$}


Let $m$ be a positive integer greater than $1$. Let $\germ g $ denote the complex Lie algebra $\germ sl _m$ of traceless $m\times m$ matrices. Let $\germ h $ denote its standard Cartan subalgebra, consisting of traceless diagonal matrices. In \cite{BBL97}, Benkart, Britten and Lemire described all the simple infinite dimensional weight $\germ g $--modules of degree $1$. Recall that a \emph{weight module} is a module for which the action of $\germ h $ is semisimple with finite multiplicities. A weight module is of degree $1$ if all its non--trivial weight spaces are $1$--dimensional.

The definition of the representations which we are interested in uses the Weyl algebra $W_m$ which is the associative algebra with generators $q_i$ and $p_i$ for $i\in \{1,\ldots , m\}$ subject to the relations $$[q_i,q_j]=0=[p_i,p_j] \mbox{ and } [p_j,q_i]=\delta_{i,j}.$$ Let $a \in \C^m$. Set 
$$P_a=\{b\in \C^m \mbox{ such that } b_i-a_i\in \Z \mbox{ for all } i\in \{1,\ldots , m\} \mbox{ and } b_i<0 \iff a_i<0\}.$$ We consider the following vector space $$W(a):=\bigoplus_{b \in P_a} \: \C x(b)$$ whose basis elements $x(b)$ are parametrized by the set $P_a$.

Now we give a structure of $W_m$--module to $W(a)$. To do so, we should think of the element $x(b)$ as a monomial $x_1^{b_1}\cdots x_m^{b_m}$ and of the $q_i$ and $p_i$ as operators of multiplication and derivation. Let $\epsilon_i$ denote the $m$--tuple whose all entries are $0$ except $i$--th entry which is $1$. We then define the (almost) natural representation of $W_m$ on $W(a)$ by 
\begin{subequations}
\label{eqW}
\begin{align}
q_i\cdot x(b) =& \left\{\begin{array}{cc}
	(b_i+1)x(b+\epsilon_i) & \mbox{if } b_i\in \Z_-\\
	x(b+\epsilon_i) & \mbox{otherwise}
	\end{array}\right.\\
p_j\cdot x(b) =& \left\{\begin{array}{cc}
	x(b-\epsilon_j) & \mbox{if } b_j\in \Z_-\\
	b_jx(b-\epsilon_j) & \mbox{otherwise}
	\end{array}\right.
\end{align}
\end{subequations}
From theorem 2.9 in \cite{BBL97}, we know that the $W_m$--module $W(a)$ is simple.

We can embed $\germ g $ (and more generally $\germ gl _m$) into $W_m$ by sending the elementary matrix $E_{i,j}$ to $q_i p_j$. We now restrict our set of parameters. Set $\cal P _a:= \{ b\in P_a\: : \:  \sum_i \: b_i=\sum_i \: a_i\}$ and $$N(a):=\bigoplus_{b\in \cal P _a} \: \C x(b).$$ A $m$--uple $b \in \cal P _a$ is called \emph{admissible} and a vector $x(b)$ associated to $b\in \cal P _a$ is called an \emph{admissible vector}. Now we have the following

\vspace{0.2cm}

\begin{theorem}[Benkart, Britten, Lemire] \cite[proposition 2.12, theorem 5.8]{BBL97} \label{thmBBL}

The vector subspace $N(a)$ of $W(a)$ is a simple weight $\germ g $--module of degree $1$. 

Conversely if $M$ is an infinite dimensional simple weight $\germ g $--module of degree $1$, then there exists $a=(a_1,\ldots, a_m)\in \C^m$, there exist two integers $k$ and $l$ with 
\begin{itemize}
\item $a_i=-1 \mbox{ for } i=1,\ldots ,k-1,$
\item $a_i\in \C\setminus \Z \mbox{ for } i=k,\ldots ,l,$
\item and $a_i=0 \mbox{ for } i=l+1,\ldots ,m,$
\end{itemize}
such that $M\cong N(a)$.
\end{theorem}

\vspace{0.2cm}

Recall that a $\germ g $--module $M$ is \emph{cuspidal} if the action on $M$ of $E_{i,j}$ is injective for all $(i,j)$ with $i\not= j$. Using the theorem \ref{thmBBL} and equations \eqref{eqW}, one shows that the cuspidal simple weight $\germ g $--modules of degree $1$ are those $N(a)$ for which $a \in (\C-\Z) ^m$. 

From now on, we suppose $m=2n$ with $n>1$. In what follows we consider only the modules of the form $$N(\underbrace{-1,\ldots ,-1}_{n-1},a_1,a_2,\underbrace{0,\ldots ,0}_{n-1}) \mbox{ with } a_1 , a_2 \in \C -\Z.$$ We denote $N_{a_1,a_2}$ this module for short. We set $a=(-1,\ldots ,-1,a_1,a_2,0,\ldots ,0)$. Remark that if $b$ is admissible then $b_i<0$ if $i<n$ and $b_j\geq 0$ if $j>n+1$. We let $\alpha_1, \ldots , \alpha_{2n-1}$ denote the standard set of simple roots for the root system $\cal R $ of $(\germ g ,\germ h )$. Then the conditions $a_1$ and $a_2\in \C-\Z$ ensure that the root vectors $X_{\pm \alpha_n}$ act injectively on $N_{a_1,a_2}$. From the action of $q_i$ and $p_j$ on $W(a)$ we derive the action of $\germ g $ on $N_{a_1,a_2}$. For example we have:
\begin{subequations}
\label{eqXi}
\begin{eqnarray}
X_{-\alpha_i} \cdot x(b)&=& \left\{\begin{array}{cc}
	(b_{i+1}+1)x(b-\epsilon_i+\epsilon_{i+1}) & \mbox{if } i<n-1\\
	x(b-\epsilon_{n-1}+\epsilon_{n}) & \mbox{if } i=n-1\\
	b_ix(b-\epsilon_i+\epsilon_{i+1}) & \mbox{if } i\geq n\\
	\end{array}\right. \label{eqnW1}\\ 
X_{\alpha_i} \cdot x(b)&=& \left\{\begin{array}{cc}
	(b_{i}+1)x(b+\epsilon_i-\epsilon_{i+1}) & \mbox{if } i<n-1\\
	b_{n}(b_{n-1}+1)x(b+\epsilon_{n-1}-\epsilon_{n}) & \mbox{if } i=n-1\\
	b_{i+1}x(b+\epsilon_i-\epsilon_{i+1}) & \mbox{if } i\geq n\\
	\end{array}\right. \label{eqnW2}\\ 
H_{\alpha_i}\cdot x(b)&=&(b_i-b_{i+1})x(b)\label{eqnW3}
\end{eqnarray}
\end{subequations}

Remark from this action that the root vectors $X_{\alpha_i}$ with $i\not= n$ act trivially on elements of the form $x(-1,\ldots ,-1,a_1+b,a_2-b,0,\ldots 0)$, with $b\in\Z$.



\section{Highest weight vectors for the action of $\germ b $ on $N_{a_1,a_2}$}


Set $\theta:=\{\alpha_1,\ldots , \alpha_{2n-1}\}\setminus \{\alpha_n\}$ and let $\langle \theta \rangle$ denote the roots in $\cal R $ which are linear combinations of the simple roots in $\theta$. To each root $\alpha \in \cal R $, we denote by $\germ g ^{\alpha}$ the associated root space. Then we associate the following Levi subalgebra $\germ l _{\theta}:=\germ h \oplus \left(\oplus_{\alpha \in \langle \theta \rangle} \: \germ g ^{\alpha}\right)$. This Levi subalgebra is the following set of matrices:
$$\germ l _{\theta} = \left\{ \left(\begin{array}{c|c}
A & 0\\
\hline
0 & B\\
\end{array}\right), \: A,B \in \mathcal{M}_n(\C), \: tr(A+B)=0\right\}.$$
Remark that the semisimple part of $\germ l _{\theta}$ is the sum of two copies of $\germ sl _n$. We denote by $\germ l _{\theta}^+$ the set of all the matrices $\left(\begin{array}{c|c}
A & 0\\
\hline
0 & B\\
\end{array}\right)$ in $\germ l _{\theta}$ with $A$ and $B$ upper triangular with zeros on the diagonal. We denote by $\germ h (\theta)$ the Cartan subalgebra of this semisimple part consisting of diagonal matrices. From this decomposition we can construct a dual pair $(\germ a ,\germ b )$ (which is $C$--admissible in the terminology of \cite{Ru94}). Define 
$$\germ b := \left\{ \left( \begin{array}{c|c}
A & 0\\
\hline
0 & A\\
 \end{array} \right), \quad \mbox{with } A \in \germ sl _n\right\} .$$ Let $\germ h _n$ denote the following Cartan subalgebra of $\germ b $:  
$$\left(\begin{array}{c|c}
D & 0\\
\hline
0 & D\\
\end{array}\right), \quad \mbox{with } D \: \mbox{diagonal and traceless}.$$ The commutant of $\germ  b $ is easily seen to be $\germ a :=\langle X,H,Y\rangle ,$ where 
$$X=\sum_{i=1}^n \: X_{-(\alpha_i+\cdots + \alpha_{n+i-1})}, \quad Y=\sum_{i=1}^n \: X_{\alpha_i+\cdots + \alpha_{n+i-1}}, \quad \mbox{and}$$ 
$$H=H_{\alpha_1}+2H_{\alpha_2}+\cdots +nH_{\alpha_n}+(n-1)H_{\alpha_{n+1}}+\cdots +2H_{\alpha_{2n-2}}+H_{\alpha_{2n-1}}.$$ The Lie algebra $\germ a $ is isomorphic to $\germ sl _2$ while $\germ b $ is isomorphic to $\germ sl _n$. It is easy to see that $(\germ a , \germ b )$ is a dual pair. Note that $\germ b \subset \germ l _{\theta}$ and that the center of $\germ l _{\theta}$ is the Cartan subalgebra $\germ h _{\theta}:=\C H$ of $\germ a $.

\vspace{0.2cm}

Now we describe the action of $\germ l _{\theta}$ on the module $N_{a_1,a_2}$.

\begin{proposition}\label{decomp:l}
As a $\germ l _{\theta}$--module, the module $N_{a_1,a_2}$ decomposes into a direct sum of infinite dimensional simple highest weight modules. Moreover the decomposition is multiplicity free and the highest weight vectors are the $x(-1,\ldots ,-1,a_1+b,a_2-b,0,\ldots ,0)$ for $b\in \Z$, whose highest weight under $\germ h (\theta)\times \germ h _{\theta}$ is $(\underbrace{0,\ldots , 0}_{n-2},-1-a_1-b,a_2-b,\underbrace{0,\ldots ,0}_{n-2})\otimes (-(n-1)+a_1-a_2+2b)$.
\end{proposition}
\dem From equation \eqref{eqnW2}, we conclude that the only vectors annihilated by $\germ l _{\theta}^+$ are exactly the admissible vectors $x(b)$ with $b=(-1,\ldots ,-1,a_1+b_n,a_2+b_{n+1},0,\ldots ,0)$ such that $b_n\in \Z$, $b_{n+1}\in \Z$, and $b_n+b_{n+1}=0$. By using equation \eqref{eqnW1} we check that for every $x(b')$ there is an element $u\in \cal U (\germ l _{\theta}^+)$ such that $u\cdot x(b')$ is a non zero multiple of some $x(b)$ with $b$ as above. Therefore the highest weight module $\cal U (\germ l _{\theta})x(b)$ corresponding to $x(b)$ is in fact simple and $N_{a_1,a_2}$ is the direct sum of the $\cal U (\germ l _{\theta})x(b)$ where $b$ is of the form $(-1,\ldots ,-1,a_1+b_n,a_2-b_n,0,\ldots ,0)$ for some integer $b_n$. The weight of this $x(b)$ is easily computed using equation \eqref{eqnW3}.

\cqfd

\begin{remark}
The modules $N_{a_1,a_2}$ give the exact list of the simple modules in the category $\cal O _{\Delta,\theta}(\germ sl _N)$ of weight--$\germ sl _{2n}$--modules $M$ which satisfy the following:
\begin{enumerate}
\item The action of $E_{n,n+1}$ and $E_{n+1,n}$ on $M$ is injective.
\item As a $\germ l _{\theta}$--module, $M$ is a direct sum of simple highest weight modules.
\end{enumerate}
See \cite{To09n} for a general definition of the categories $\cal O _{\Delta,\theta}$, and their properties.
\end{remark}

\vspace{0.2cm}

For later use, let us compute the action of the root vectors occuring in $X \mbox{ and } Y \in \germ a $:

\begin{lemma}\label{lemA}
For $i \in \{1,\ldots , n\}$, we have 
\begin{subequations}
\label{eqXYA}
\begin{eqnarray}
X_{\alpha_i+\ldots +\alpha_{n+i-1}}x(b)&=&\left\{\begin{array}{cc}
(b_i+1)b_{n+i}x(b-\epsilon_{n+i}+\epsilon_i) & \mbox{ if } i<n,\\
b_{2n}x(b-\epsilon_{2n}+\epsilon_n) & \mbox{ if } i=n,\\
\end{array}\right.\\
X_{-(\alpha_i+\ldots +\alpha_{n+i-1})}x(b)&=&\left\{\begin{array}{cc}
x(b+\epsilon_{n+i}-\epsilon_i) & \mbox{ if } i<n,\\
b_{n}x(b+\epsilon_{2n}-\epsilon_n) & \mbox{ if } i=n,\\
\end{array}\right.
\end{eqnarray}
\end{subequations}
\end{lemma}
\dem The root vector $X_{\alpha_i+\ldots +\alpha_{n+i-1}}$ correspond to the elementary matrix $E_{i,n+i}$. The root vector $X_{-(\alpha_i+\ldots +\alpha_{n+i-1})}$ correspond to the elementary matrix $E_{n+i,i}$. Now the lemma follows from equations \eqref{eqW}. 

\cqfd

\vspace{0.2cm}

The first step toward a correspondence is to understand the action of $\germ b $ on $N_{a_1,a_2}$. We set $\germ b ^+$ the subalgebra of $\germ b $ consisting of the matrices $\left( \begin{array}{c|c}
A & 0\\
\hline
0 & A\\
 \end{array} \right)$ in $\germ b $ such that $A$ is upper triangular with zero diagonal. From proposition \ref{decomp:l}, we obtain that the action of $\germ b ^+$ is locally finite. We will now investigate the subspace 
$$M_0:=N_{a_1,a_2}^{\germ b ^+}=\{x \in N_{a_1,a_2} \: : \:  X\cdot x=0, \: \forall \: X \in \germ b ^+\}.$$ Define $X_i:=X_{\alpha_i}+X_{\alpha_{n+i}}$, $X_{-i}:=X_{-\alpha_i}+X_{-\alpha_{n+i}}$, and $H_i=[X_i,X_{-i}]$. From equation \eqref{eqnW2} we get 
\begin{equation}\label{eqnB}
X_i\cdot x(b)  =  \left\{ \begin{array}{cc}
(b_i+1)x(b+\epsilon_i-\epsilon_{i+1}) + b_{n+i+1}x(b+\epsilon_{n+i}-\epsilon_{n+i+1}) & \mbox{if } i<n-1\\
b_n(b_{n-1}+1)x(b-\epsilon_n+\epsilon_{n-1})+b_{2n}x(b+\epsilon_{2n-1}-\epsilon_{2n}) & \mbox{if } i=n-1\\
\end{array} \right.
\end{equation}

\begin{lemma}\label{lemB1}
Let $x=\sum_{k} \: \lambda_k x(b^k) \in M_0$ be a weight vector for $\germ h _{\theta}\oplus \germ h _n$. Then there exist indices $k_0$ and $k_1$  such that 
$$x(b^{k_0})=x(-1,\ldots ,-1,b_n^{k_0},\ldots , b_{2n}^{k_0})$$ and
$$x(b^{k_1})=x(b_1^{k_1},\ldots , b_{n+1}^{k_1},0,\ldots , 0).$$
\end{lemma}
\dem If $b$ is admissible, then $b_1<0$. Let $i_1$ be an index such that $\lambda_{i_1}\not=0$ and $b_1^{i_1}$ is maximal among the possible values of the different $b_1^k$ occuring in $x$. Suppose $b_1^{i_1}\not=-1$. Then applying $X_1$ to $x(b^{i_1})$ gives according to equation (\ref{eqnB}) the following sum of two vectors:
$$(b^{i_1}_1+1)x(b^{i_1}+\epsilon_1-\epsilon_2)+b^{i_1}_{n+2}x(b^{i_1}+\epsilon_{n+1}-\epsilon_{n+2}).$$
The first summand is a vector $x(b')$ such that $b'_1=b_1^{i_1}+1>b_1^{i_1}$. The second summand is of the form $x(b'')$ with $b''_1=b^{i_1}_1$. But by our hypothesis on $x$, we have $X_1\cdot x=0$. We look at the other occurences of the vector $x(b')$ in $X_1\cdot x$. By the maximality of $b_1^{i_1}$, this vector only occur as the first summand of $X_1 \cdot x(b^{i_1})$. Thus its coefficient in $X_1\cdot x$ is $b^{i_1}_1+1$ which is non zero by our hypothesis on $b_1^{i_1}$, contradicting the fact that $x\in M_0$. Hence $b_1^{i_1}=-1$.

If $b$ is admissible, then we also have $b_2<0$. So let now $i_2$ denote an index such that $\lambda_{i_2}\not=0$, $b_1^{i_2}=-1$ and $b_2^{i_2}$ is maximal among the possible values of the different $b_2^k$ occuring in $x$ and subject to the condition that $b_1^k=-1$. We apply the same reasonning using equation \eqref{eqnB} for the action of $X_2$ to prove that $b_2^{i_2}=-1$. Applying then $X_3,\ldots X_{n-1}$ we get an index $k_0$ satisfying the condition of the lemma, i.e. $b_1^{k_0}=\cdots =b_{n-1}^{k_0}=-1$.

If $b$ is admissible, then $b_{2n}\geq 0$. Therefore to find $k_1$ we do the same thing starting from the action of $X_{n-1}$ to a vector $x(b^{j_1})$ such that $\lambda_{j_1}\not=0$ and $b_{2n}^{j_1}$ is minimal among the possible $b_{2n}^{k}$. We prove that necessarily $b_{2n}^{j_1}=0$. Applying successively $X_{n-2},\ldots , X_1$ we obtain an index $k_1$ satisfying the condition of the lemma, i.e. $b_{2n}^{k_1}=\cdots =b_{n+1}^{k_1}=0$.

\cqfd

\begin{corollary}\label{corB2}
Let the notations be as in lemma \ref{lemB1}. Then there are integers $c_n, c_{n+1}, c_{2n}$ with $c_{2n}\geq 0$ such that
$$x(b^{k_0})=x(-1,\ldots , -1,a_1+c_n,a_2+c_{n+1},0,\ldots , 0, c_{2n})$$ and 
$$x(b^{k_1})=x(-1-c'_1,-1, \ldots , -1,a_1+c'_n,a_2+c'_{n+1},0,\ldots , 0),$$ with $c_n+c_{n+1}+c_{2n}=0$, $c'_1=c_{2n}$, $c'_n=c_n+c_{2n}$ and $c'_{n+1}=-c_n$.
\end{corollary}
\dem Set $b_n^{k_0}=a_1+c_n$, $b_{n+1}^{k_0}=a_2+c_{n+1}$, $b_{n+i}^{k_0}=c_{n+i}$ for $2\leq i\leq n$. Set also $b_n^{k_1}=a_1+c'_n$, $b_{n+1}^{k_1}=a_2+c'_{n+1}$, $b_i^{k_1}=-1-c'_i$ for $1\leq i\leq n-1$. Therefore, 
$$b^{k_0}=(-1,\ldots ,-1,a_1+c_n,a_2+c_{n+1},c_{n+2},\ldots ,c_ {2n}),$$
$$b^{k_1}=(-1-c'_{1},\ldots , -1-c'_{n-1},a_1+c'_n,a_2+c'_{n+1},0,\ldots ,0).$$
The two vectors $x(b^{k_0})$ and $x(b^{k_1})$ should be admissible vectors and should have the same weight with respect to the Cartan subalgebras of $\germ a $ (since it commutes with $\germ b $) and to the Cartan subalgebra of $\germ b $ (generated by $H_1, \ldots ,H_{n-1}$). This gives rise to the following equations
\begin{align*}
\mbox{admissibility of } x(b^{k_0}): \quad c_n+c_{n+1}+\ldots + c_{2n} & = 0 \\
\mbox{admissibility of } x(b^{k_1}): \quad -(c'_1+\ldots + c'_{n-1})+c'_n+c'_{n+1} & = 0 \\
\mbox{$\germ h _{\theta}$--weight of } x(b^{k_0}) \mbox{ and } x(b^{k_1}): \quad c_n-(c_{n+1}+\ldots +c_{2n}) & = -(c'_1+\ldots + c'_{n-1})+c'_n-c'_{n+1} \\
\begin{array}{r}
\mbox{weight of } x(b^{k_0}) \mbox{ and } x(b^{k_1})\\
 \mbox{ under } H_1,\ldots ,H_{n-1} 
\end{array} : \quad \left\{ \begin{array}{r}
c_{n+1}-c_{n+2}\\
c_{n+2}-c_{n+3}\\
 \\
c_{2n-2}-c_{2n-1}\\
-c_n+c_{2n-1}-c_{2n}
\end{array} \right. & \begin{array}{c}
=\\
=\\
\vdots \\
=\\
=
\end{array} \begin{array}{c}
c'_2-c'_1+c'_{n+1}\\
c'_3-c'_2\\
 \\
c'_{n-1}-c'_{n-2}\\
-c'_{n-1}-c'_n
\end{array}  
\end{align*}
Some calculations show that the unique solution in the $c'_i$'s variables of this system is 
$$\left\{ \begin{array}{ccc}
c'_1&=&-c_n-c_{n+1}\\
c'_i&=&-c_{n+i}, \: 1<i<n\\
c'_n&=&c_{2n}+c_n\\
c'_{n+1}&=&-c_n\\
\end{array} \right.$$
But now as the vectors $x(b^{k_0})$ and $x(b^{k_1})$ must be admissible, we should have $c_{n+i}\geq 0$ for $2\leq i\leq n$ and $c'_j\geq 0$ for $1\leq j\leq n-1$. This imposes that $c_{n+i}=c'_i=0$ for $2\leq i\leq n-1$. Then taking into account that $c_n+c_{n+1}+\cdots +c_{2n}=0$ again by admissibity of $b^{k_0}$, we finally obtain
\begin{equation}
c_{n+1}=-c_n-c_{2n}, \quad c_{2n}\geq 0, \quad c'_1=c_{2n}, \quad c'_{n}=c_n+c_{2n}, \quad c'_{n+1}=-c_n
\end{equation}

\cqfd

\vspace{0.2cm}

From now on, we shall write $x(b^{k_0})=x(-1,\ldots ,-1,a_1+b,a_2-b-c,0,\ldots ,0,c)$ where $b$ and $c$ are integers and $c\geq 0$. Let us investigate which admissible vectors $x(b)$ have the same weight with respect to $\germ h _{\theta}\oplus \germ h _n$ than $x(b^{k_0})$. This is the following lemma:

\begin{lemma}\label{lemB3}
Set $x:=x(-1-b_1,\ldots , -1-b_{n-1}, a_1+b_n,a_2+b_{n+1},b_{n+2},\ldots , b_{2n})$ an admissible vector. Then $x$ has the same weight than $x(b^{k_0})$ under the action of $\germ h _{\theta}\oplus \germ h _n$ if and only if 
$$\left\{ \begin{array}{ccc}
b_{n+i} & = & b_{i}, \quad 2\leq i \leq n-1\\
b_n & = & b+b_1+\cdots +b_{n-1}\\
b_{n+1} & = & b_1-(b+c)\\
b_{2n} & = & c-(b_1+\cdots +b_{n-1})
\end{array} \right.$$
\end{lemma}
\dem As in the proof of corollary \ref{corB2}, we write down the equations obtained by expressing the admissibility of $x$ and the fact that $x$ and $x(b^{k_0})$ have the same weight under $\germ h _{\theta}$ and $\germ h _n$:
\begin{align}
-(b_1+\cdots +b_{n-1})+b_n+\cdots +b_{2n} & = 0 \label{eqLB3a}\\
-(b_1+\ldots + b_{n-1})+b_n-(b_{n-1}+\cdots +b_{2n}) & = 2b \label{eqLB3b}\\
\left\{ \begin{array}{r}
b_2-b_1+b_{n+1}-b_{n+2}  \\
b_3-b_2+b_{n+2}-b_{n+3}  \\
  \\
b_{n-1}-b_{n-2}+b_{2n-2}-b_{2n-1} \\
-b_{n-1}-b_n+b_{2n-1}-b_{2n}   
\end{array} \right.
\begin{array}{c}
=\\
=\\
\vdots  \\
=\\
=
\end{array} & \begin{array}{c}
-b-c\\
0\\
\\
0\\
-b-c
\end{array}\label{eqLB3c}
\end{align}
Then we set $\tilde{b}_i=b_i-b_{n+i}$ for $1\leq i \leq n-1$ and $\tilde{b}_n=b_n+b_{2n}$. We rewrite equations \eqref{eqLB3a} and \eqref{eqLB3c} in the new variables $\tilde{b}_i$:
\begin{align*}
-(\tilde{b}_1+\cdots +\tilde{b}_{n-1})+\tilde{b}_n & = 0\\
\left\{ \begin{array}{r}
\tilde{b}_2-\tilde{b}_1\\
\tilde{b}_3-\tilde{b}_2\\
 \\
\tilde{b}_{n-1}-\tilde{b}_{n-2}\\
-\tilde{b}_{n-1}-\tilde{b}_n
\end{array} \right. \begin{array}{c}
=\\
=\\
\vdots \\
=\\
=
\end{array} \begin{array}{c}
-b-c\\
0\\
 \\
0\\
-b-c
\end{array} & 
\end{align*}
The unique solution of this system in the $\tilde{b}_i$'s variables is $\tilde{b}_2=\cdots =\tilde{b}_{n-1}=0$, $\tilde{b}_1=b+c$, $\tilde{b}_n=b+c$. Therefore, we have
$$b_n+b_{2n}=c+b, \: b_{n+1}=b_1-c-b, \mbox{ and } b_{n+i}=b_i \mbox{ for } 2\leq i\leq n-1.$$ Then using equation \eqref{eqLB3b}, we express $b_{n+i}$ for $i\geq 0$ in the $b_j$'s variables for $1\leq j \leq n-1$, which gives the lemma.

\cqfd

\begin{corollary}\label{corB4}
Let $x\in M_0$ be a weight vector with respect to $\germ h _{\theta}\oplus \germ h _n$. Then there are two integer $b$ and $c$ such that $c\geq 0$ and 
$$x=\sum_{k_i\geq 0, \: |\underline{k}| \leq c} \: \lambda_{\underline{k}} \: x(-1-k_1,\ldots , -1-k_{n-1}, a_1+b+|\underline{k}|,a_2-b-c+k_1,k_2,\ldots , k_{n-1},c-|\underline{k}|),$$ where $\underline{k}=(k_1,\ldots , k_{n-1})\in \N^{n-1}$, $|\underline{k}|=\sum_i \: k_i$ and $\lambda_{\underline{k}}\in \C$.

If $n>2$, its $\germ h _n$--weight is $(a_2-b-c,0,\ldots , 0,-1-a_1-b-c)$ and its $\germ h _{\theta}$--weight is $a_1-a_2+2b-(n-1)$. If $n=2$ then its $\germ h _n$--weight is $(-1-a_1+a_2-2(b+c))$ and its $\germ h _{\theta}$--weight is $a_1-a_2+2b-1$.
\end{corollary}
\dem From lemma \ref{lemB1} and corollary \ref{corB2} we know that $x=\sum_i \: \lambda_i \: x(b^i)$ and that there is an index $i_0$ such that $b^{i_0}=(-1, \ldots,-1,a_1+b,a_2-b-c,0,\ldots, 0,c)$ for some integers $b$ and $c$ with $c\geq 0$. Then the lemma \ref{lemB3} asserts that the others $x(b^i)$ occuring in $x$ are of the form $x(b^i)=x(-1-k_1,\ldots ,-1-k_{n-1},a_1+b+(k_1+\cdots +k_{n-1}),a_2+k_1-b-c,k_2,\ldots ,k_{n-1},c-(k_1+\cdots +k_{n-1}))$. These vectors should also be admissible. Therefore we must have
$$k_i \in \N, \quad \mbox{and } c-(k_1+\cdots +k_{n-1})\geq 0.$$ This is the corollary.

\cqfd

\vspace{0.2cm}

\begin{proposition}\label{propB5}
Let $x$ be as in the corollary \ref{corB4}. Write 
$$x=\sum_{k_i\geq 0, \: |\underline{k}| \leq c} \: \lambda_{\underline{k}} \: x(-1-k_1,\ldots , -1-k_{n-1}, a_1+b+|\underline{k}|,a_2-b-c+k_1,k_2,\ldots , k_{n-1},c-|\underline{k}|),$$ for some integers $b$ and $c$ with $c\geq 0$. Then 
$\lambda_{\underline{k}}=\kappa (\underline{k})\lambda_{\underline{0}},$ where $\lambda_{\underline{0}}\in \C$ and 
$$
\kappa (\underline{k})=\binom{k_1+k_2}{k_1}\cdots \binom{k_1+\cdots+k_{n-1}}{k_1+\cdots +k_{n-2}}\frac{\prod_{j=1}^{k_1+\cdots +k_{n-1}} (c+1-j)}{(k_1+\cdots +k_{n-1})!\prod_{j=1}^{k_1+\cdots +k_{n-1}} (a_1+b+j)}.$$
Conversely, if $$x=\sum_{k_i\geq 0, \: |\underline{k}| \leq c} \: \lambda_{\underline{k}} \: x(-1-k_1,\ldots , -1-k_{n-1}, a_1+b+|\underline{k}|,a_2-b-c+k_1,k_2,\ldots , k_{n-1},c-|\underline{k}|)$$ with $\lambda_{\underline{k}}=\kappa (\underline{k})\lambda_{\underline{0}}$,  then $x\in M_0$.
\end{proposition}
\dem From equation \eqref{eqnB} we have 
$$X_1\cdot x = $$
$$\sum_{\underline{k}}\: \lambda_{\underline{k}}\Big{[}-k_1x(-k_1,-2-k_2,-1-k_3,\ldots ,-1-k_{n-1},a_1+b+|\underline{k}|,a_2-b-c+k_1,k_2,\ldots , k_{n-1},c-|\underline{k}|)$$
$$ + k_2x(-1-k_1,\ldots ,-1-k_{n-1},a_1+b+|\underline{k}|,a_2-b-c+k_1+1,k_2-1,k_3,\ldots ,k_{n-1},c-|\underline{k}|)\Big{]}.$$

Let $\underline{k}\in\{(k_i)_{1\leq i\leq n-1} \: : \: k_i\geq 0 \mbox{ and } \sum_i \: k_i \leq c\}.$ Suppose $k_1>0$. Let $\underline{k}'=(k'_i)_{1\leq i\leq n-1}$ be such that $k'_1=k_1-1$, $k'_2=k_2+1$, and $k'_i=k_i$ otherwise. We look at the coefficient of 
$$x(-k_1,-2-k_2,-1-k_3,\ldots ,-1-k_{n-1},a_1+b+|\underline{k}|,a_2-b-c+k_1,k_2,\ldots , k_{n-1},c-|\underline{k}|)$$ in the expression of $X_1\cdot x$. We find $$-k_1\lambda_{\underline{k}}+(k_2+1)\lambda_{\underline{k}'}.$$ As $x\in M_0$, we have $X_1\cdot x=0$. Therefore we should have $-k_1\lambda_{\underline{k}}+(k_2+1)\lambda_{\underline{k}'}=0$, i.e. $\lambda_{\underline{k}}=\frac{k_2+1}{k_1}\lambda_{\underline{k}'}$. By induction we find that $\lambda_{\underline{k}}=\binom{k_1+k_2}{k_1}\lambda_{\underline{k}^1}$ where $\underline{k}^1=(0,k_1+k_2,k_3,\ldots ,k_{n-1})$. We then look at the coefficient of 
$$x(-1,-k_2-k_1,-2-k_3,-1-k_4,\ldots ,-1-k_{n-1},a_1+b+|\underline{k}|,a_2-b-c+k_1,k_2,\ldots , k_{n-1},c-|\underline{k}|)$$ in $X_2\cdot x(\underline{k}^1)$. This allows us to express $\lambda_{\underline{k}^1}$ from $\lambda_{\underline{k}^2}$ where $\underline{k}^2=(0,0,k_1+k_2+k_3,k_4\ldots ,k_{n-1})$. More precisely, we get $\lambda_{\underline{k}^1}=\binom{k_1+k_2+k_3}{k_1+k_2}\lambda_{\underline{k}^2}$. Then using successively the action of $X_3,\ldots ,X_{n-1}$ on $x$, we express $\lambda_{\underline{k}}$ from $\lambda_{\underline{0}}$.

The converse is easy.

\cqfd

From now on, we denote by $x(b,c)$ the vector in $M_0$ obtained in proposition \ref{propB5} such that $\lambda_{\underline{0}}=1$. We also denote by 
$$x_{\underline{k}}(b,c)=x(-1-k_1,\ldots , -1-k_{n-1}, a_1+b+|\underline{k}|,a_2-b-c+k_1,k_2,\ldots , k_{n-1},c-|\underline{k}|)$$ for $\underline{k}=(k_1,\ldots ,k_{n-1})$ such that $k_i\in \N$ and $|\underline{k}|=k_1+\cdots +k_{n-1}\leq c$.

\begin{corollary}\label{corB6}
For $2\leq i\leq n-2$ we have $X_{-i}\cdot x(b,c)=0$.
\end{corollary}
\dem From equation \eqref{eqnW2}, we see that $X_{-i}$ acts trivially on $x_{\underline{k}}(b,c)$ if and only if $k_i=k_{i+1}=0$. If the action is non trivial, we have $X_{-i}\cdot  x_{\underline{k}}(b,c)=-(k_{i+1})x_{\underline{k}'}(b,c)+k_ix_{\underline{k}''}(b,c)$ where 
$$k'_j=\left\{ \begin{array}{cc}
k_j & \mbox{if } j\not=i \mbox{ or } i+1\\
k_i-1 & \mbox{if } j=i\\
k_{i+1}+1 & \mbox{if } j=i+1\\
\end{array}\right. ,$$ and 
$$k''_j=\left\{\begin{array}{cc}
k_j & \mbox{if } j\not=n+i \mbox{ or } n+i+1\\
k_{n+i}-1 & \mbox{if } j=n+i\\
k_{n+i+1}+1 & \mbox{if } j=n+i+1
\end{array}\right. .$$ We now look at the occurences of $x_{\underline{k}'}(b,c)$ in $X_{-i}\cdot x(b,c)$. It appears in the expression of $X_{-i}\cdot  x_{\underline{k}}(b,c)$ as we already mentionned and in the second summand of the expression of $X_{-i}\cdot x_{\underline{l}}(b,c)$ where 
$$l_j=\left\{\begin{array}{cc}
k_j & \mbox{if } j\not=n+i \mbox{ or } n+i+1\\
k_{n+i}+1 & \mbox{if } j=n+i\\
k_{n+i+1}-1 & \mbox{if } j=n+i+1
\end{array}\right. .$$ Therefore the coefficient of $x_{\underline{k}'}(b,c)$ in $X_{-i}\cdot x(b,c)$ is $(-k_{i+1})\lambda_{\underline{k}}+l_i\lambda_{\underline{l}}$. Using the expression of $\lambda_{\underline{k}}$ given by proposition \ref{propB5}, we see that this coefficient is $0$.

\cqfd



\section{A Howe type correspondence for $N_{a_1,a_2}$}


\subsection{Generic case}

Now to state and prove a Howe--type correspondence for the module $N:=N_{a_1,a_2}$ we need to compute the action of $\germ a $ on $M_0$.

\begin{lemma}\label{lemA2}
Let $b$ and $c$ be two integers with $c\geq 0$. Then 
\begin{enumerate}
\item $X\cdot x(b,c)$ is a non zero element of $M_0$ which is equal to a multiple of $x(b-1,c+1)$.
\item The vector $Y\cdot x(b,c)$ is non zero if and only if $c(a_1-a_2+2b+c-(n-2))\not=0$. In this case, it is equal to a multiple of $x(b+1,c-1)$.
\end{enumerate}
\end{lemma}
\dem As the action of $\germ a $ commutes with the action of $\germ b $, we get that $Y\cdot x(b,c)$ and $X\cdot x(b,c)$ both are in $M_0$. Let us set $b'=b+1$ and $c'=c-1$. Then using lemma \ref{lemA}, we compute $Y\cdot x(b,c)$. We obtain 
$$
Y\cdot x(b,c)=
$$
$$
\sum \: \lambda_{\underline{k}}(a_2-b'-c'+k_1)(-k_1)x(-k_1,-1-k_2,\ldots , -1-k_{n-1},a_1+b'+|\underline{k}|-1,$$
$$a_2-b'-c'+k_1-1,k_2,\ldots , k_{n-1},c'-|\underline{k}|+1)$$
$$
+ \sum_{i=2}^{n-2} \: \sum \: \lambda_{\underline{k}}(-k_i^2)x(-1-k_1,-1-k_2,\ldots ,-1-k_{i-1},-k_i,-1-k_{i+1},\ldots , -1-k_{n-1},a_1+b'+|\underline{k}|-1,$$
$$a_2-b'-c'+k_1,k_2,\ldots , k_{i-1}, k_i-1,k_{i+1},\ldots , k_{n-1},c'-|\underline{k}|+1)$$
$$
 + \sum \: \lambda_{\underline{k}}(c'+1-|\underline{k}|)x(-1-k_1,\ldots , -1-k_{n-1},a_1+b'+|\underline{k}|,a_2-b'-c'+k_1,k_2,\ldots , k_{n-1},c'-|\underline{k}|)
$$
As $Y\cdot x(b,c) \in M_0$, it should be a linear combination of some $x(b'',c'')$. But each $x(b'',c'')$ contains a vector of the form $x_{\underline{0}}(b'',c'')=x(-1,\ldots ,-1,a_1+b'',a_2-b''-c'',0,\ldots ,0,c'')$. The only such vector in the expression of $Y\cdot x(b,c)$ is $x_{\underline{0}}(b',c')$. Therefore if $Y\cdot x(b,c)$ is non zero then it is a multiple of $x(b',c')$. To see when it is zero, it is enough to compute the coefficient of $x_{\underline{0}}(b',c')$ in the above equation. The first sum gives a contribution equal to $-(a_2-b-c+1)\lambda_{\underline{\epsilon_1}}$, the second sum gives $\sum_{i=2}^{n-2} \: -\lambda_{\underline{\epsilon_i}}$ and the last sum gives $c$ (recall that the coefficient $\lambda_{\underline{0}}$ of $x(b,c)$ was set equal to $1$). Using proposition \ref{propB5} which allows us to express the $\lambda_{\epsilon_i}$'s, we find the global contribution:
$$\frac{c}{a_1+b+1}(a_1-a_2+2b+c-(n-2)).$$ Thus $Y\cdot x(b,c)\not=0$ if and only if $c(a_1-a_2+2b+c-(n-2))\not=0$.
 
We apply the same method for $X\cdot x(b,c)$. The coefficient of $x_{\underline{0}}(b-1,c+1)$ in the expression of $X\cdot x(b,c)$ is $a_1+b$, which is non zero since $a_1\in \C\setminus \Z$. This gives the lemma.
 
\cqfd

\begin{corollary}\label{corB7}
Assume $a_1-a_2\not\in \Z$. Then for all integers $b$ and $c$ with $c\geq 0$ the $\germ b $--module generated by $x(b,c)$ is a simple  highest weight module.
\end{corollary}
\dem The $\germ b $--module $\cal U (\germ b )x(b,c)$ generated by $x(b,c)$ is a highest weight module and is therefore indecomposable. Thus it is simple if and only if $x(b,c)$ is the only highest weight vector in $\cal U (\germ b )x(b,c)$, up to a scalar multiple. Another highest weight vector in $\cal U (\germ b )x(b,c)$ would be of the form $x(b',c')$. But then the vectors $x(b,c)$ and $x(b',c')$ would have the same $\germ h _{\theta}$--weight. Using the $\germ h _{\theta}$--weight given in corollary \ref{corB4}, we see that necessarily $b'=b$. Thus if $x(b,c')\in \cal U (\germ b )x(b,c)$, then there is an element $u\in \cal U (\germ b ^-)$ such that $u\cdot x(b,c)=x(b,c')$. This implies first that $c\leq c'$ since the vectors $X_{-i}$ can only increase the $2n$--th component of every admissible vector. Hence if $c\not= c'$, we have $c<c'$. Then from the hypothesis $a_1-a_2\not\in \Z$ and from lemma \ref{lemA2}, we get $Y^{c+1}\cdot x(b,c)=0$ and $Y^{c+1}\cdot x(b,c')\not=0$. As $Y\in \germ a $ commutes with $u\in \cal U (\germ b )$, we should have $uY^{c+1}x(b,c)=Y^{c+1}x(b,c')$. This is a contradiction.

\cqfd

\vspace{0.2cm}

Let $L(\lambda)$ denote the simple highest weight $\germ sl _2$--module with highest weight $\lambda$ and $L(a_2-b,0,\ldots , 0,-1-a_1-b)$ denote the simple highest weight $\germ sl _n$--module with highest weight $(a_2-b,0,\ldots , 0,-1-a_1-b)$.

\begin{theorem}\label{thmAB1}
Assume $a_1-a_2\not\in \Z$. Then we have the following decomposition of $N_{a_1,a_2}$ as a $\germ b \oplus \germ a $--module:
\begin{enumerate}
\item If $n=2$, $$N_{a_1,a_2}=\bigoplus_{b\in \Z} \: L(-1-a_1+a_2-2b)\otimes L(a_1-a_2+2b-1).$$
\item If $n>2$, $$N_{a_1,a_2}=\bigoplus_{b\in \Z} \: L(a_2-b,0,\ldots , 0,-1-a_1-b)\otimes L(a_1-a_2+2b-(n-1)).$$
\end{enumerate}
\end{theorem}
\dem We know from proposition \ref{decomp:l} that the module $N_{a_1,a_2}$ is $\germ l _{\theta}^+$--locally finite. Thus it is also $\germ b ^+$--locally finite. Therefore for every vector $v$ in $N_{a_1,a_2}$ there is an element $u$ of $\cal U (\germ b ^+)$ such that $u\cdot v$ is in $M_0$. From corollary \ref{corB7}, we know that each weight vector in $M_0$ spans a simple highest weight $\germ b $--module. Thus $N_{a_1,a_2}=\oplus_{b\in \Z, c\in \N} \: \cal U (\germ b )x(b,c)$ is the decomposition of $N_{a_1,a_2}$ into simple $\germ b $--modules. The hypothesis $a_1-a_2\not\in \Z$ ensures that $Y\cdot x(b,c)=0$ if and only if $c=0$. Then we get the following chain of $\germ b $--modules :
$$0\rightharpoonup \cal U (\germ b )x(b,0) \rightleftharpoons U (\germ b )x(b-1,1) \rightleftharpoons \cdots \rightleftharpoons U (\germ b )x(b-k,k) \rightleftharpoons \cdots ,$$ where $\rightharpoonup$ stands for the action of $X$ and $\leftharpoondown$ for the action of $Y$.
Thanks to lemma \ref{lemA2}, we conclude that this chain is a simple $\germ a \oplus \germ b $--module which is then by corollary \ref{corB4} isomorphic to $L(-1-a_1+a_2-2b)\otimes L(a_1-a_2+2b-1)$ if $n=2$ and to $L(a_2-b,0,\ldots , 0,-1-a_1-b)\otimes L(a_1-a_2+2b-(n-1))$ if $n>2$.

\cqfd

\vspace{0.2cm}

As a consequence of this theorem, we find the following Howe--type correspondence in the $''$generic$''$ case $a_1-a_2 \not\in\Z$, namely:
$$L(-1-a_1+a_2-2b)\leftrightarrow L(a_1-a_2+2b-1), \quad \mbox{if } n=2,$$
$$L(a_2-b,0,\ldots , 0,-1-a_1-b)\leftrightarrow L(a_1-a_2+2b-(n-1)), \quad \mbox{if } n>2.$$

\subsection{Non--generic case}

Let us now consider the non--generic case $a_1-a_2\in \Z$. From lemma \ref{lemA2}, we find that $Y\cdot x(b,c)=0$ if and only if $c=0$ or $a_1-a_2+2b+c=n-2$ and that $X\cdot x(b,c)$ is always non zero. We also get that the $\germ a $--module generated by $x(b,0)$ is a highest weight module of highest weight $a_1-a_2+2b-(n-1)$. As a vector space, this module is $\bigoplus_{k\in \N} \: \C x(b-k,k)$. The vector $x(b-k,k)$ for $k>0$ is annihilated by $Y\in \germ a $ if and only if $a_1-a_2+2b-k=n-2$. Therefore there is at most one $k$ for which $x(b-k,k)$ is a highest weight vector for $\germ a $. Thus we have shown the following:

\begin{corollary}\label{corA3}
Assume $a_1-a_2\in \Z$. Let $b\in \Z$.
\begin{enumerate}
\item If $a_1-a_2+2b-(n-2)\leq 0$ then the $\germ a $--module generated by $x(b,0)$ is irreducible.
\item If $a_1-a_2+2b-(n-2)>0$, then the $\germ a $--module generated by $x(b,0)$ has length $2$:
$$\cal U (\germ a )x(b,0)\supset \cal U (\germ a )x(b-(a_1-a_2+2b-(n-2)),a_1-a_2+2b-(n-2))\supset \{0\},$$ where $\cal U (\germ a )x(b-(a_1-a_2+2b-(n-2)),a_1-a_2+2b-(n-2))$ is a simple highest weight $\germ a $--module (of weight $a_1-a_2+2b-(n-1)$) and the quotient 
$$\cal U (\germ a )x(b,0)/\cal U (\germ a )x(b-(a_1-a_2+2b-(n-2)),a_1-a_2+2b-(n-2))$$ is a simple highest weight $\germ a $--module of weight $a_1-a_2+2b-(n-1)$. In this case the $\germ a $--module $\cal U (\germ a )x(b,0)$ is isomorphic to the Verma module $V(a_1-a_2+2b-(n-1))$ of highest weight $a_1-a_2+2b-(n-1)$.
\end{enumerate}
\end{corollary}

\vspace{0.2cm}

We can now use the same method as in corollary \ref{corB7} and prove the following result:

\begin{corollary}\label{corB8}
Assume $a_1-a_2\in \Z$. Let $b\in \Z$ such that $a_1-a_2+2b-(n-2)\geq 0$. Then for all $c\in \N$, the $\germ b $--module generated by $x(b,c)$ is a simple highest weight $\germ b $--module of weight $(a_2-b-c,0,\cdots ,0,-1-a_1-b-c)$
\end{corollary}
\dem The proof of corollary \ref{corB7} can also be applied in this case because the hypothesis on $b$ together with corollary \ref{corA3} ensures that $Y\cdot x(b,c)$ is non zero as soon as $c>0$.

\cqfd

\vspace{0.2cm}

In general, the same argument shows that the only vector $x(b,c')$ that can belong to the $\germ b $--module generated by $x(b,c)$ satisfies $c'>c$ and $a_1-a_2+2b+c+c'=n-2$ (see the proof of corollary \ref{corB7}). Now we prove the following

\begin{proposition}\label{propB9}
Let $b\in \Z$ such that $a_1-a_2+2b-(n-2)< 0$. Let $c\in \N$ be such that $a_1-a_2+2b+c-(n-2)=0$. Then the $\germ b $--module $\cal U (\germ b )x(b,c)$ is contained in the $\germ b $--module $\cal U (\germ b )x(b,0)$. Moreover, the latter has length $2$, with the following composition serie:
$$\cal U (\germ b )x(b,0) \supset \cal U (\germ b )x(b,c)\supset \{0\}.$$
\end{proposition}
\dem To prove the first assertion, it suffices to find $u\in \cal U (\germ b )$ such that $u\cdot x(b,0)=\alpha x(b,c)$, with $\alpha \in \C ^*$. We define the following element. Let $Z'\in \germ b $ be the element corresponding to the matrix $E_{n,1}+E_{2n,n+1}$ and let $$Z''=\sum_{i=1}^{n-2} \left(E_{i+1,1}+E_{n+i+1,n+1}\right)\left(E_{n,i+1}+E_{2n,n+i+1}\right)\in\cal U (\germ b ).$$ Remark that $[Z',Z'']=0$. For $\lambda \in \C$, we denote by $Z_{\lambda}$ the vector $Z'+\lambda Z''$. Then $[Z_{\lambda},Z_{\lambda'}]=0$ for all $\lambda, \lambda'\in \C$. We list now some computations of brackets :
\begin{align*}
[X_1,Z_{\lambda}] = & \left(E_{n,2}+E_{2n,n+2}\right)\left(\lambda (n-2)-1+\lambda H_1\right)\\
 & -\lambda \sum_{i=2}^{n-2} \: \left(E_{n,i+1}+E_{2n,n+i+1}\right)\left(E_{i+1,2}+E_{n+i+1,n+2}\right),\\
[X_i,Z_{\lambda}] = & 0, \quad \mbox{for } 2\leq i \leq n-2,\\
[X_{n-1},Z_{\lambda}] = & \left(E_{n-1,1}+E_{2n-1,n+1}\right)\left(1+\lambda H_{n-1}\right)\\
 & +\lambda \sum_{i=1}^{n-3}\: \left(E_{i+1,1}+E_{n+i+1,n+1}\right)\left(E_{n-1,i+1}+E_{2n-1,n+i+1}\right),
\end{align*}
\begin{align*}
[E_{n,2}+E_{2n,n+2},Z_{\lambda}] = & \lambda Z'\left(E_{n,2}+E_{2n,n+2}\right),\\
[H_1,Z_{\lambda}] = & -Z_{\lambda},\\
[E_{i+1,2}+E_{n+i+1,n+2},Z_{\lambda}] = & 0,
\end{align*}
\begin{align*}
[E_{n-1,1}+E_{2n-1,n+1},Z_{\lambda}] = & -\lambda Z'\left(E_{n-1,1}+E_{2n-1,n+1}\right),\\
[H_{n-1},Z_{\lambda}]= & -Z_{\lambda},\\
[E_{n-1,i+1}+E_{2n-1,n+i+1},Z_{\lambda}] = & 0.
\end{align*}

For $1\leq i\leq c$, set $\lambda _i=\frac{1}{a_1+b+i}$ (note that $a_1+b+i$ is non zero since $a_1\in \C\setminus \Z$). Set also $Z_i=Z_{\lambda_i}$ and $Z=Z_1\cdots Z_c$. Now set $x=Z\cdot x(b,0)$. We show that $x$ is a highest weight vector for $\germ b $. We already know that $x(b,0)$ is a highest weight vector for $\germ b $. Thus we have to show that $ad(X_i)(Z)\cdot x(b,0)=0$ for $1\leq i \leq n-1$. From the relations above, we already find that $ad(X_i)(Z)=0$ for $2\leq i \leq n-2$. Let us compute $ad(X_1)(Z)\cdot x(b,0)$. We obtain 
$$
ad(X_1)(Z)\cdot x(b,0)=\: [X_1,Z_1]Z_2\cdots Z_c\cdot x(b,0) + \cdots + Z_1\cdots Z_{c-1}[X_1,Z_c]\cdot x(b,0).$$ In the expression of $[X_1,Z_i]$, appear the vectors $\left(E_{k+1,2}+E_{n+k+1,n+2}\right)$. But we have seen that these vectors commute with all the $Z_j$'s. Moreover from corollary \ref{corB6}, we get that $\left(E_{k+1,2}+E_{n+k+1,n+2}\right)$ acts trivially on $x(b,0)$. Therefore the only part in the expression of $[X_1,Z_i]$ that can give a non trivial contribution in the expression of $ad(X_1)(Z)\cdot x(b,0)$ is $\left(E_{n,2}+E_{2n,n+2}\right)\left(\lambda_i (n-2)-1+\lambda_i H_1\right)$. Thus,
\begin{align*}
ad(X_1)(Z)\cdot x(b,0)= & \left(E_{n,2}+E_{2n,n+2}\right)\left(\lambda_1 (n-2)-1+\lambda_1 H_1\right)Z_2\cdots Z_c\cdot x(b,0)+\cdots \\
 & +Z_1\cdots Z_{c-1}\left(E_{n,2}+E_{2n,n+2}\right)\left(\lambda_c (n-2)-1+\lambda_c H_1\right)\cdot x(b,0).
\end{align*}
From our previous computations we also get that 
$$H_1Z_k\cdots Z_c=-(c+1-k)Z_k\cdots Z_c+Z_k\cdots Z_cH_1.$$ Thus we have 
\begin{align*}
ad(X_1)(Z)\cdot x(b,0) = & Z_2\cdots Z_c\left(\lambda_1 (n-2)-1+\lambda_1 H_1-(c-1)\lambda_1\right)\cdot x(b,0)+\cdots \\
 & +Z_1\cdots Z_{c-1}\left(\lambda_c (n-2)-1+\lambda_c H_1\right)\cdot x(b,0).
\end{align*}
Now from corollary \ref{corB4}, we get that $H_1\cdot x(b,0)=(a_2-b)x(b,0)$. Then using the definition of $c$, we conclude that $\left(\lambda_k (n-2)-1+\lambda_k H_1-(c-k)\lambda_k\right)\cdot x(b,0)=0$, which in turn expresses that $ad(X_1)(Z)\cdot x(b,0)=0$.

Now since $[Z_{\lambda},Z_{\lambda'}]=0$, we have also that $Z=Z_c\cdots Z_1$. Then we compute $ad(X_{n-1})(Z)\cdot x(b,0)$ using the same method as above and prove that $ad(X_{n-1})(Z)\cdot x(b,0)=0$. Therefore we have proved that $Z\cdot x(b,0)$ is a highest weight vector for $\germ b $ (note that $Z\cdot x(b,0)$ is a weight vector because $Z$ is a weight vector in $\cal U (\germ b )$).

It only remains to show that $Z\cdot x(b,0)\not=0$. To do so, we compute the coefficient of $x_{\underline{0}}(b,c)$ in the expression of $Z\cdot x(b,0)$. We have seen that $Z'$ commute with $Z''$. So $Z$ is a homogeneous polynomial of degree $c$ in the two variables $Z'$ and $Z''$. Remark that $Z''$ cannot increase the $2n$--th component of the admissible vectors. Therefore the only monomial in the expression of $Z$ that can give the admissible vector $x_{\underline{0}}(b,c)$ when acting on $x(b,0)$ is $Z'^c$, whose coefficient in the polynomial $Z$ is $1$. After some computations we find that the coefficient of $x_{\underline{0}}(b,c)$ in $Z\cdot x(b,0)$ is $(a_2-b)(a_2-b-1)\cdots (a_2-b-(c-1))$. This is non zero since $a_2 \in \C\setminus \Z$. Therefore from proposition \ref{propB5}, we conclude that $Z\cdot x(b,0)$ is a non zero multiple of $x(b,c)$.

From our choice of $c$ and lemma \ref{lemA2}, we show that $Y\cdot x(b,c+k)\not=0$ for $k>0$. Then we can apply the same proof as in \ref{corB7} to show that $\cal U (\germ b )x(b,c)$ is a simple $\germ b $--module. Thus the $\germ b $--module $\cal U (\germ b )x(b,0)$ is a highest weight module (containing the simple $\germ b $--module $\cal U (\germ b )x(b,c)$), which has therefore a composition serie of finite length consisting of highest weight modules. But we have remarked that $\cal U (\germ b )x(b,0)$ cannot contain any other highest weight vector than the linear combinations of $x(b,0)$ and $x(b,c)$ (see above the statement of this proposition). From this we conclude that $$\cal U (\germ b )x(b,0) \supset \cal U (\germ b )x(b,c)\supset \{0\}$$ is the composition serie of $\cal U (\germ b )x(b,0)$.

\cqfd

\vspace{0.2cm}

We can now state and prove the following Howe--type correspondence in the non--generic case:

\begin{theorem}\label{thmAB2}
Assume $a_1-a_2\in \Z$ and $n>2$. Let $b\in \Z$. Then we have the following correspondence:
\begin{itemize}
\item If $a_1-a_2+2b-(n-2)=0$, then $\cal U (\germ b )x(b,0)$ is a simple $\germ b $--module isomorphic to $L(a_2-b,0,\ldots ,0,-1-a_1-b)$ and we have $L(a_2-b,0,\ldots ,0,-1-a_1-b) \leftrightarrow L(-1)$.
\item If $a_1-a_2+2b-(n-2)>0$, then $\cal U (\germ b )x(b,0)$ is a simple $\germ b $--module isomorphic to $L(a_2-b,0,\ldots ,0,-1-a_1-b)$ and we have $L(a_2-b,0,\ldots ,0,-1-a_1-b) \leftrightarrow V(a_1-a_2+2b-(n-1))$.
\item If $a_1-a_2+2b-(n-2)<0$, then $\cal U (\germ b )x(b,0)$ is an indecomposable $\germ b $--module of length $2$ and we have $\cal U (\germ b )x(b,0)\leftrightarrow L(a_1-a_2+2b-(n-1))$.
\end{itemize}
\end{theorem}
\dem This is a consequence of corollaries \ref{corA3} and \ref{corB8} and of proposition \ref{propB9}.

\cqfd

\begin{theorem}\label{thmAB3}
Assume $a_1-a_2\in \Z$ and $n=2$. Let $b\in \Z$. Then we have the following correspondence:
\begin{itemize}
\item If $a_1-a_2+2b=0$, then $\cal U (\germ b )x(b,0)$ is a simple $\germ b $--module isomorphic to $L(-1)$ and we have $L(-1) \leftrightarrow L(-1)$.
\item If $a_1-a_2+2b>0$, then $\cal U (\germ b )x(b,0)$ is a simple $\germ b $--module isomorphic to $L(a_2-a_1-2b-1)$ and we have $L(a_2-a_1-2b-1) \leftrightarrow V(a_1-a_2+2b-1)$.
\item If $a_1-a_2+2b<0$, then $\cal U (\germ b )x(b,0)$ is an indecomposable $\germ b $--module of length $2$ isomorphic to $V(a_2-a_1-2b-1)$ and we have $V(a_2-a_1-2b-1) \leftrightarrow L(a_1-a_2+2b-1)$.
\end{itemize}
\end{theorem}
\dem This is a consequence of corollaries \ref{corA3} and \ref{corB8} and of proposition \ref{propB9}.

\cqfd

Remark that in both case we have a correspondence of the infinitesimal caracteres. Let us give an interpretation of the theorem \ref{thmAB2}. First we have the following decomposition of $N_{a_1,a_2}$ as a $\germ b $--module:
$$N_{a_1,a_2}=\sum_{b\in \Z, c\in \N}\: \cal U (\germ b )x(b,c).$$ But we have seen that each $\cal U (\germ b )x(b,c)$ is either simple or has length $2$. Let us then consider the semisimplification $N_{a_1,a_2}^s$ of $N_{a_1,a_2}$ obtained by changing those $\cal U (\germ b )x(b,c)$ which are indecomposable by their composition factors. The space $N_{a_1,a_2}^s$ is still a $\germ b $--module (but not a $\germ g $--module anymore) and we have the following branching rules:
$$N_{a_1,a_2}^s=\bigoplus_{b\in \Z, c\in \N}\: L(a_2-b-c,0,\ldots ,0,-1-a_1-b-c).$$ But it is clear that $N_{a_1,a_2}^s$ still caries an action of $\germ a $ induced by the action of $\germ a $ on $N_{a_1,a_2}$. From theorem \ref{thmAB2}, we find the following Howe--type correspondence for this $\germ b \oplus \germ a $--module $N_{a_1,a_2}^s$:
\begin{itemize}
\item If $a_1-a_2+2b-(n-2)=0$, then we have $L(a_2-b,0,\ldots ,0,-1-a_1-b) \leftrightarrow L(-1)$.
\item If $a_1-a_2+2b-(n-2)>0$, then we have $L(a_2-b,0,\ldots ,0,-1-a_1-b) \leftrightarrow V(a_1-a_2+2b-(n-1))$.
\item If $a_1-a_2+2b-(n-2)<0$, then we have $L(a_2-b,0,\ldots ,0,-1-a_1-b) \leftrightarrow L(a_1-a_2+2b-(n-1))$.
\end{itemize}
Note that in this correspondence some $\germ a $--modules are not simple modules. So we can consider the corresponding semisimplification $N_{a_1,a_2}^{ss}$ of $N_{a_1,a_2}^s$ (which we call the bi--semisimplification of $N_{a_1,a_2}$). In this module, we get the following correspondence:
\begin{itemize}
\item If $a_1-a_2+2b-(n-2)=0$, then we have $L(a_2-b,0,\ldots ,0,-1-a_1-b) \leftrightarrow L(-1)$.
\item If $a_1-a_2+2b-(n-2)>0$, then we have $L(a_2-b,0,\ldots ,0,-1-a_1-b) \leftrightarrow L(a_1-a_2+2b-(n-1))\oplus L(-(a_1-a_2+2b-(n-3))$.
\item If $a_1-a_2+2b-(n-2)<0$, then we have $L(a_2-b,0,\ldots ,0,-1-a_1-b) \leftrightarrow L(a_1-a_2+2b-(n-1))$.
\end{itemize}
Note that this is no more a one to one correspondence. We can also give an interpretation of theorem \ref{thmAB3} in the same spirit. The final correspondence in the $''$bi--semisimplification$''$ of $N_{a_1,a_2}$ is in this case the following:
\begin{itemize}
\item If $a_1-a_2+2b=0$, then we have $L(-1) \leftrightarrow L(-1)$.
\item If $a_1-a_2+2b>0$, then we have $L(a_2-a_1-2b-1) \leftrightarrow L(a_1-a_2+2b-1)\oplus L(-(a_1-a_2+2b+1))$.
\item If $a_1-a_2+2b<0$, then we have $L(a_2-a_1-2b-1) \leftrightarrow L(a_1-a_2+2b-1)$.
\end{itemize}


\bibliographystyle{plain}
\bibliography{biblio}

\end{document}